\documentclass[10pt,twoside]{article}
\usepackage{amsmath,amsfonts,amssymb,amscd,amsthm,amsbsy,epsf}
\usepackage{Latex-document}
\newtheorem{thm}{Theorem}[section]

\def\C{{\mathcal C}}
\def\F{{\mathcal F}}
\def\K{{\mathcal K}}

\def\osc{\operatorname{osc}}
\def\ep{\varepsilon}
\def\negint{\mathop{\int\mkern-19mu{-}}\nolimits}
\renewcommand{\thesection}{\arabic{section}}

\def\volumeno{I}

\title{\bf Non Linear Elliptic Theory and \vskip -2mm the Monge-Ampere Equation\vskip 6mm}
\author{Luis A. Caffarelli\thanks{Department of Mathematics, University of Texas at Austin,
Austin, TX 78712, USA. E-mail: caffarel@math.utexas.edu}
\vspace*{-0.5cm}}
\date{\vspace{-8mm}}

\begin{document}

\maketitle

\thispagestyle{first} \setcounter{page}{179}

\begin{abstract}

\vskip 3mm

The Monge-Ampere equation, plays a central role in the theory of
fully non linear equations. In fact we will like to show how the
Monge-Ampere equation, links in some way the ideas comming from
the calculus of variations and those of the theory of fully non
linear equations.

\vskip 4.5mm


%\noindent {\bf Keywords and Phrases:}
\end{abstract}

\vskip 12mm

When learning complex analysis, it was a remarkable fact that the real
part $u$ of an analytic function, just because it satisfies the equation:
$$u_{xx} + u_{yy} =  \Delta\, u = 0 $$
(Laplace's equation)
is real analytic, and furthermore, the oscillation of $u$ in any given
domain $U$, controls {\em all\/} the derivatives of $u$, of {\em any\/}
order, in any subset $\bar U$, compactly contained in $U$.

One can give three, essentially different explanations of this phenomena.

\noindent {\bf a) Integral representations} (Cauchy integral, for
instance){\bf .} This gives rise to many of the modern aspects of
real and harmonic analysis: fundamental solutions, singular
integrals, pseudo-differential operators, etc. For our discussion,
an important consequence of this theory are the Schauder and
Calderon-Zygmund estimates.

Heuristically, they say that if we have a solution of an equation
$$A_{ij} (x) D_{ij} u = 0$$
and $A_{ij}(x)$ is, in a given functional space, a small perturbation
of the Laplacian then $D_{ij}u$ is actually in the same functional space as
$A_{ij}$.
For instance, if $[A_{ij}]$ is H\"older continuous ($C^{\alpha}(\bar U)$) and positive definite,
we can transform it to the identity (the Laplacian) at any given
point $x_0$ by an affine transformation, and will remain close to it in a neighborhood. Thus $D_{ij}u$ will also be $C^{\alpha}(\bar U)$.

\noindent {\bf b) Energy considerations.} Harmonic functions, $u$,
are also local minimizers of the Dirichlet integral
$$E(v) = \int(\nabla v)^2\,dx\ .$$
That is, if we change $u$ to $w$, in $\bar U \subset\subset U$
$$E(w)|_{\bar U} \ge E(u)|_{\bar U}\ .$$
This gives rise to the theory of calculus of variations
(minimal surface, harmonic maps, elasticity, fluid dynamics).

One is mainly concerned, there, with equations (or systems) of the form
\begin{equation}\label{eq1}
D_i F_i (\nabla u,X) = 0 \ .
\end{equation}
For instance, in the case in which $u$ is a local minimizer of
$$E(u) = \int \F (\nabla u,X)\,dx$$
\eqref{eq1} is simply the Euler-Lagrange equation associated to $E$:
$$F_i = \nabla_p \F\ .$$
If we attempt to write \eqref{eq1} in second derivatives form, we get
$$F_{i,j} (\nabla u,X) D_{ij} u + \cdots = 0\ .$$

This strongly suggests that in order for the variational problem to be
``elliptic'', like the Laplacian, $F_{i,j}$ should be positive definite,
that is $\F$ should be strictly convex.

 It also leads to the natural strategy of showing that $\nabla u$, that
   in principle is only in $L^2$ (finite energy), is in
fact  H\"older continuous. Reaching this regularity
allows us  to apply  the (linear) Schauder theory.

That  implies  $D_{ij} u$ is $C^\alpha(\bar U)$, thus $\nabla u$ is
$C^{1,\alpha}(\bar U)$, and so on (the bootstrapping method).

The difficulty with this approach is that solutions, $u$, are
invariant under $R^{n+1}$-dialations of their graphs.

This fact keeps the class of Lipschitz functions (bounded gradients) invariant.
There is no reason, thus, to expect that this equation will ``improve''
under dialations. The fact that $\nabla u$ is indeed H\"older continuous
is the celebrated De Giorgi's theorem, that solved the nineteenth Hilbert's problem:

De Giorgi looked at the equation that first derivatives, $u_\alpha$ satisfy
$$D_i F_{ij} (\nabla u) D_j u_\alpha = 0 .$$

He thought of $F_{ij} (\nabla u)$ as elliptic coefficients $A_{ij} (x)$ that had no regularity whatsoever, and he proved that any solution $w$ of
$$D_i A_{ij} (x) D_j w = 0 $$
was  H\"older continuous
$$\|w\|_{C^\alpha (\bar U)} \le C\| w\|_{L^2 (U)}\ .$$

De Giorgi's theorem is in fact a linear one, but for a new invariant class
of equations. No matter how the solution (and the equation) is renormalized,
it stays far from the constant coefficient theory, and a radically new idea
surfaces:  if we have a class of functions for which at every scale, in
some average sense, the function controls its derivatives
(the energy inequality), further regularity follows.

%The passage from $\nabla u\in L^2(U)$ to $\nabla u\in C^\alpha(\bar U)$, is, of course,
%no trivial matter.
%It is the celebrated De Giorgi theorem, that then evolved into the
%De Giorgi-Nash-Moser theory.
%In fact, the De Giorgi theorem is much more powerful than that.
%It considers a variational solution of the {\em linear\/} equation
%$$D_i A_{ij} (x) D_j w = 0 $$
%but without assuming any regularity on the coefficients $A_{ij}(x)$, only
%ellipticity, and it proves that such a $w$ is H\"older continuous.%

%Furthermore
%$$\|w\|_{C^\alpha (\bar U)} \le C\| w\|_{L^2 (U)}\ .$$
%In doing so, De Giorgi makes a {\em jump\/} of invariance classes.
%
%Previous to his work, the theory could consider only  equations of the form
%$$D_i A_{ij} (x) D_j u = 0$$
%that were a small perturbation of the Laplacian, that is, that under
%dialations become asymptotically the Laplacian.
%  with just bounded coefficients.
%We are now confronted with an equation that no matter how much we
%dilate, remains in the same class.

Finally, the third approach is

\noindent {\bf c) Comparison principle.} Two solutions $u_1,u_2$
of $\Delta u=0$ cannot ``touch without crossing''. That is, if
$u_1-u_2$ is positive it cannot become zero in some interior
point, $X_0$, of $U$.

Again, heuristically, this is because the function
$$F(D^2 u) = \Delta u = \text{Trace}[D^2u]$$
is a monotone function of the Hessian matrix $[D_{ij}u] $
and, thus, in some sense, we must have $F(D^2u_1)$ ``$>$'' $F(D^2u_2)$
at $X_0$ (or nearby).

The natural family of equations to consider in this context, is then
$$F(D^2u) =0$$
for $F$ a strictly monotone function of $D^2u$.

Such type of equations appear in differential geometry.
For instance, the coefficients of the characteristic polynomial of the
Hessian
$$P(\lambda) = \det (D^2 u - \lambda I)$$
are such equations if we restrict $D^2u$ to stay in the appropriate
set of $R^{n\times n}$.
If $\lambda_i$ denote the eigenvalues of $D^2u$
\begin{equation*}
\begin{split}
C_1 & = \Delta u = \sum \lambda_i \qquad \text{(Laplace)}\cr
C_2 & = \sum_{i\ne j} \lambda_k \lambda_j \ldots\cr
C_n & = \prod \lambda_i = \det D^2 u\qquad \text{(Monge-Ampere)}\ .
\end{split}
\end{equation*}
In the case of $C_n = \det D^2 u = \prod\lambda_i$ is a monotone function
of the Hessian provided that all $\lambda_i$'s are positive.
That is, provided that the function, $u$, under consideration is convex.

If $F(D^2 u,X) $ is {\em uniformly elliptic\/}, that is, if $F$ is strictly
monotone as a function of the Hessian, or in differential form,
$$F_{ij} (M) = D_{m_{ij}} F$$
is uniformly positive definite, then solutions of $F(D^2 u)$ are
$C^{1,\alpha}(\bar U)$.
As in the divergence case, this is because first derivatives $u_\alpha$ satisfy an
elliptic operator,
$$F_{ij}(D^2 u) D_{ij} u_\alpha =0$$
now in non divergence form. As long as we do not have further information on
$D^2 u$, we must think  again of $F_{ij}$ as
 bounded measurable coefficients.

The De Giorgi type theorem for $a_{ij}(x) D_{ij} u_\alpha =0$
is due to Krylov and Safanov, and states again that solutions of such
an equation are H\"older continuous.

We point out that, again this result has ``jumped'' invariance classes.
Rescaling of $a_{ij}(x)$ does not improve them.
Unfortunately, this is not enough to ``bootstrap'', as in the divergence
case:
The coefficients, $A_{ij}(x) = F_{ij} (D^2 u)$, depend on second
derivatives.
If we will manage to prove that $D^2u$ is H\"older continuous, then, from
equation \eqref{eq1}, $D_\alpha u$ would be  $C^{2,\alpha}(\bar U)$, i.e.,
$u$ would be $C^{3,\alpha}(\bar U)$ and we could improve and improve.

To prove this, once more convexity reappears.
If $F(D^2u)$ is concave (or convex) then all pure second derivatives are
sub (or super) solutions of the linearized operator.
This, together with the fact that $D^2u$ lies in the surface
$F(D^2u)$, implies the H\"older continuity of $D^2u$, and, by the
bootstrapping argument $u$ is as smooth as $F$ allows.

\subsection*{The Monge-Ampere equation and optimal transportation}

\vskip-5mm \hspace{5mm}

We would like now to turn our attention to the Monge-Ampere
equation
\begin{equation*}
\det D^2 u = \prod \lambda_i = f(x,u,\nabla u)\ .
\end{equation*}
As pointed out before, the equation fits  in the context of elliptic
equations provided that we consider convex solutions.
That is, provided that $f$ is positive.
Further $\log\det D^2u = \sum \log \lambda_i$ is concave as function
of the $\lambda_i$ and thus is a concave function of $D^2u$.
Unfortunately  $\det D^2u$ {\em is not\/} uniformly
strictly convex.

For instance if we prescribe
\begin{equation*}
\det D^u  = \prod\lambda_i = 1
\end{equation*}
ellipticity deteriorates as one of the $\lambda$'s goes to infinity
and some other is forced to go to zero.
This difficulty is compensated by two fundamental facts.
\begin{itemize}
\item[1)] The rich family of invariances that the Monge-Ampere equation
enjoys.
\item[2)] Its ``hidden'' divergence structure.
\end{itemize}

The divergence structure is due to the fact that $\det D^2u$ can be
thought of as the Jacobian of the gradient map: $X\to \nabla u$.
Thus for any domain $\bar U$
\begin{equation*}
\int_{\bar U} \det D^2 u \,dx = \text{Vol}(\nabla u(\bar U)).
\end{equation*}
But if $\bar U \subset\subset U$, $u$ being convex implies that
\begin{equation*}
(\nabla u) |_{\bar U} \le C\osc u|_U.
\end{equation*}
This gives us a sort of ``energy inequality''  that controls a
positive quantity of $D^2 u$ by the oscillation of $u$:
\begin{equation*}
\int_{\bar U} \det D^2 u \le C(\bar U,U) (\osc u)^n.
\end{equation*}

\subsection*{Invariances}

\vskip-5mm \hspace{5mm}

The Monge-Ampere equation is invariant of course, under the the
standard families of transformations:
\begin{itemize}
\item[a)] Rigid motions, $R$:
\begin{equation*}
\det D^2 u (Rx) = f(Rx),
\end{equation*}
\item[b)] Translations:
\begin{equation*}
\det D^2 u (x + v) = f( x+v),
\end{equation*}
\item[c)] Quadratic dialations:
\begin{equation*}
\det D^2 \frac1{t^2} u(tx) = f(tx).
\end{equation*}
But also
\item[d)] Monge-Ampere is invariant under any affine transformation $A$,
of determinant one:
\begin{equation*}
\det D^2 u(Ax) = f(Ax)\ .
\end{equation*}
\end{itemize}

If $f$ is, for instance, in one of the following classes:
\begin{itemize}
\item[a)] $f$ constant,
\item[b)] $f$ close to constant ($|f-1| \le \ep$),
\item[c)] $f$  bounded away from zero and infinity
($0< \frac1{\sigma} \le f\le \sigma$),
\end{itemize}
any of the transformations above gives a new $u$ in the {\em same class\/}
of solutions.

For instance, if $u$ is a solution of
\begin{equation*}
\det D^2 u = 1
\end{equation*}
then, $u(\ep x,\frac1{\ep}y)$ is also a solution of the {\em same  equation}.
But this has dramatically ``deformed'' the graph of $u$.
It is then almost unavoidable that there are singular solutions
(Pogorelov).

In fact, for $n\ge 3$, one can construct convex solutions $u$ that contain
a line their graph and are not differentiable in the direction
transversal to that line, solutions of
\begin{equation*}
\det D^2 u = f(x)
\end{equation*}
with $f$ a smooth positive function.

Fortunately, this geometry can only be inherited from the boundary
of the domain.

\begin{thm}\label{thm:main}
If in the domain $U\subset R^n$
\begin{itemize}
\item[a)] $\frac1{\sigma} \le \det D^2 u \le \sigma$,
\item[b)] $u\ge 0$,
\item[c)] The set $\Gamma = \{u=0\}$ is not a point, then $\Gamma$ is
generated as ``convex combinations'' of its boundary points
\begin{equation*}
\Gamma = \text{ convex envelope of }\ \Gamma\cap \partial U\ .
\end{equation*}
\end{itemize}
\end{thm}

A corollary of this theorem is that
\begin{itemize}
\item[a)] If we can ``cut a slice'' of the graph of $u$, with a hyperplane
 $l(x)$ so that  the support $S$ of $(u-l)^-$
is compactly contained in $U$, then
$u$ is, inside $S$, both $C^{1,\alpha}$ regular and also $C^{1,\alpha}$-
strictly convex, i.e., separates from any of its supporting planes with
polynomial growth.
\end{itemize}

This is the equivalent of De Giorgi's and Krylov-Safanov result
(remember that the $C^\alpha$ theorems were applied to the {\em
derivatives\/} of the solutions of the non-linear equations under
consideration).

Note that by an affine transformation and a dilation we can always
renormalize the  support of the ``slice''  $S$
 to be equivalent to the unit ball of $R^n$:
$B_1 \subset S\subset B_n$.

After this normalization, it is possible to reproduce for $u$ all the
classical estimates we had for the Laplacian:
\begin{itemize}
\item[a)] (Calderon-Zygmund).
If $f$ is close to constant $(|f-1|<\ep)$, then $D^2u \in L^p(B_{1/2})$
($p=p(\ep)$ goes to infinity when $\ep$ goes to zero).
\item[b)]
If $f\in C^{k,\alpha}$ (has up to $k$ derivatives H\"older continuous) then
$u\in C^{k+2,\alpha}$ (all second derivatives of $u$ are $C^{k,\alpha}$.
\end{itemize}

Note that $f$ plays, for Monge-Ampere, simultaneously the role of
``right hand side'' and ``coefficients'' due to the structure of
its non-linearity.

\subsection*{The Monge-Ampere equation and optimal transportation
(the Monge problem)}

\vskip-5mm \hspace{5mm}

The Monge-Ampere equation has many applications, not only in
geometry, but also in applied areas: optimal design of antenna
arrays, vision, statistical mechanics, front formation in
meteorology, financial mathematics.

Many of these applications are related to optimal transportation and the
Wasserstein metric between probability distributions.
In the discrete case, optimal transportation consists of the following.

We are given two sets of $k$ points in $R^n$: $X_1,\ldots,X_k$ and
$Y_1,\ldots, Y_k$, and want to map the $X$'s onto the $Y$'s, i.e., we
look at all one-to-one functions $Y(X_j)$.
But we want to do so, minimizing some transportation costs
\begin{equation*}
\C = \sum_j C\Big( Y(X_j) - X_j\Big)\ .
\end{equation*}

For our discussion $C(X-Y) = \frac12 |X-Y|^2$.
It is easy to see that the minimizing map must be the gradient
(subdifferential) of a convex potential $\varphi$.

In the continuous case, instead of having $k$-points we have two
probability densities, $f(X) \,dX$ and $g(Y)\,dY$
and we want to consider those (admissible) maps $Y(X)$ that
``push forward'' $f$ to $g$.

Heuristically that means that in the change of variable formula, we can
substitute
\begin{equation*}
g(Y(X)) \det D_X Y(X) \mbox{``=''} f(X).
\end{equation*}
A weak formulation, substitutes the map $Y(X)$, by a  joint
probability density $\nu (X,Y)$ with marginals $f(X)\,dX $ and
$g(Y) \,dY$, i.e.,
\begin{equation*}
\begin{split}
f(X_0) = \int d_Y \nu (X_0,Y),\cr g(Y_0) = \int d_X\nu (X,Y_0).
\end{split}
\end{equation*}
(We don't ask the ``map'' to be one-to-one any more, the image of $X_0$
may now spread among ``many $Y$'s''.

Among all such $\nu$, we want to maximize correlation
\begin{equation*}
\K = \int\langle X,Y\rangle d\nu (X,Y)
\end{equation*}
or minimize cost
\begin{equation*}
\C = \int \frac12 |X-Y|^2 \,d\nu (X,Y),
\end{equation*}
$\sqrt{\C}$ defines a metric, the Wasserstein metric among probability
densities.

Under mild hypothesis, we have the

\begin{thm}\label{thm:mild}
The unique optimal $\nu_0$ concentrates in a graph (is actually a
one-to-one map, $Y(X)$). Further $Y(X)$ is the subdifferential of
a convex potential $\varphi$, i.e., $Y(X) = \nabla \varphi$.
Heuristically, then, $\varphi$ must satisfy the Monge-Ampere
equation
\begin{equation*}
g(\nabla \varphi) \det D^2 \varphi = f(X).
\end{equation*}
\end{thm}

For several reasons, the weak theory does not apply in general,
but one can still prove, for instance:

\begin{thm}\label{thm:weak}
If $f$ and $g$ never vanish or if the supports of $f$ and $g$ are convex sets,
the map $Y(X)$ is ``one derivative better'' than $f$ and $g$.
\end{thm}

\subsection*{Some applications and current issues}

\vskip-5mm \hspace{5mm}

{\bf a)}
 It was pointed out by Otto, that the Wasserstein metric can be used to
describe the evolution of several of the classical ``diffusion'' equations:
heat equation, porous media, lubrication.

The idea is that a diffusion process for one equation with conservation of
mass, consists of the balance of two factors:
trying to minimize distance between consecutive distributions
($u(x,t_k)$ and $u(x,t_{k+1})$), plus trying to flatten or smooth
(diffuse), $u(x,t_{k+1})$.

This fact has allowed to prove rates of decay to equilibrium in
many of the classical equations, as well as a number of new
phenomena. The fine relations between the discrete and continuous
problems is an evolving issue (rate of convergence, regularity of
the discrete problems, etc.).

{\bf b)}
 Another family of problems, coming both from geometry and optimal
transportation concerns the study of several issues on solutions
of Monge-Ampere equations in periodic or random media.

{\bf b$_1$)   Liouville type theorems:} We start with a theorem of
Calabi of Liouville type: Given a global convex solution of
Monge-Ampere equation, $\det D^2u=1$, $u$ must be a quadratic
polynomial. Suppose now that instead of RHS equal to one, we have
a general RHS, $f(x)$. Given a global solution, to discover its
behavior at infinity we may try to ``shrink it'' through quadratic
transformations:
\begin{equation*}
u_\ep = \ep^2 u \left( \frac{x}{\ep}\right)\ ,\ \text{ satisfies }\
\det D^2 u_\ep = \left( \frac{x}{\ep}\right)\ .
\end{equation*}
Suppose now that $f$ averages out at infinity, for instance $f$ is
periodic. Then due to the ``divergence structure'' of Monge-Ampere
$u_\ep$ should converge to a quadratic polynomial.

\begin{thm}\label{thm:Liouville1}
Given a RHS\ $f(x)$, periodic, with average $\negint f=a$
\begin{itemize}
\item[i)] Given any quadratic polynomial $P$ with $\det D^2P=a$, there
exists a unique periodic function $w$, such that
\begin{equation*}
\det D^2 (P+w) = f(x)
\end{equation*}
($w$ is a ``corrector'' in homogenization language).
\item[ii)] Conversely (Liouville type theorem):
Given a global solution $u$, it must be of the form $P+w$.
\end{itemize}
\end{thm}

What are the implications for homogenization?
What can we say if $f(X,u,\nabla u)$ is periodic in $X$ and $u$?
What can we say if $f_\omega(x)$ is random in $X$?

{\bf b$_2$)  Vorticity transport:}
(2 dimensions)
Again in the periodic  context we seek a ``vorticity density'',
$\rho (X,t)$ periodic in $X$. At each time $t$, $\rho$ generates a periodic
``stream function'', $\psi (X,t)$ by the equation
\begin{equation*}
\det (I+D^2\psi) = \rho\ .
\end{equation*}
In turn, $\psi$ generates a periodic velocity field $v= - (\psi_y,\psi_x)$
that transports $\rho$:
\begin{equation*}
\rho_t + \text{div} (v\rho) =0\ .
\end{equation*}
Given some initial data $\rho_0 (x)$, what can we say about $\rho$?

If $\rho_0$ is a vorticity patch, $\rho_0 (x) = 1+\chi_\Omega$, does
it stay that way?

If we choose $\rho_0$, $\psi_0$ so that $\rho_0 = F(\psi_0)$, that
is $\det I +D^2\psi_0 = F(\psi_0)$, we have a stationary vorticity
array, i.e., $\rho (X,t) \equiv \rho_0$.

What can we say, in parallel to the classic theory of rotating fluids,
or plasma, where $\det$ is substituted by $\Delta \psi$?

{\bf c)}
Another area of research relates to optimal transportation as a natural
``map'' between probability densities.
It has been shown that optimal transportation explains naturally interpolation
properties of densities (of Brunn Minkowski type), monotonicity properties
(like correlation inequalities that express in which way the probability
density, $g$, is shifted in some cone of directions with respect to $f$),
and concentration properties of $g$ versus $f$ (in which sense for instance,
a log concave perturbation of a Gaussian is more concentrated than a
Gaussian).

Of particular interest would be to understand optimal transportation as
dimension goes to infinity. Since convex potentials are very stable objects,
this would provide, under some circumstances, an ``infinite dimensional''
 change of variables formula between probability densities.

{\bf d)}
Finally, one of my favorite problems is to understand the geometry
of optimal transportation in the case in which the cost function
$C(X-Y)$ is still strictly convex, but not quadratic.
In that case, the optimal map is still related to a potential that satisfies
\begin{equation*}
\det (I+D (F_j (\nabla \psi))) = \cdots
\end{equation*}
where $F_j$ is now the gradient of the convex conjugate to $C$.

At this point, we have come full circle and we are now in a higher hierarchy,
in a sort of  Lagrangian
version of the Euler-Lagrange equation from the calculus of variations.

In fact if we put an epsilon in front of $D$ and linearize,
\begin{equation*}
\det (I+\ep D(F_j (\nabla \psi))) = 1 + \ep \text{ Trace} (D(F_j (\nabla \psi)))+ O(\ep^2)
= 1+\ep \text{ div } F_j (\nabla \psi)\ + O(\ep^2).
\end{equation*}

Bibliographical references can be found in the books of J.
Gilbarg-N. Trudinger, L.C. Evans and L.A. Caffarelli-X. Cabre for
nonlinear PDE's; T. Aubin, I. Bakelman and C. Gutierrez for the
Monge-Ampere equation, and the recent surveys by L. Ambrosio and
C. Villani for optimal transportation.

\label{lastpage}

\end{document}